    \documentclass[reqno]{amsproc}
    \usepackage{amsfonts}
    \usepackage{amssymb}
   \usepackage{amsmath}
   \usepackage{euscript}
   \usepackage{amsfonts}
\usepackage{amssymb}

\usepackage{euscript}
\usepackage{amsmath}

\usepackage{txfonts}

   \DeclareMathOperator{\okr}{{\stackrel{{\scriptscriptstyle{def}}}{=}}}
   \DeclareMathOperator{\D}{d\!} \DeclareMathOperator{\E}{e}
   \DeclareMathOperator{\I}{i} \DeclareMathOperator{\lin}{lin}
   \DeclareMathOperator{\diag}{{diag}}

    \def\zzz{\mathbb Z}

    \newcommand{\bbinom}[2]{\genfrac{[}{]}{0pt}{}{#1}{#2}}
   \def\Le{\leqslant}
   \def\Ge{\geqslant}
   \def\ddc{\mathcal D}
    \def\dz#1{\mathcal D({#1})}
    \def\gw{^*}
    \def\hhc{\mathcal H}
    \def\is#1#2{\langle#1,#2\rangle}
    \def\jd#1{\mathcal N(#1)}
    \def\kkc{\mathcal K}
   \def\Ld#1{\mathcal L^2(#1)}
    \def\liczp#1{{${#1}^{\text {\rm o}}$}}
    \def\naw#1{{\rm(}{#1}{\rm )}}
    
    \def\nul{\{0\}}
    \def\ob#1{\mathcal R(#1)}
   \def\od{^{-1}}
    
    \def\res#1{|_{#1}}
    \def\rres#1{\!\!\upharpoonright_{#1}}
    \def\wyn{\text{$\implies$ }}
    \def\zb#1#2{\{{#1};\;{#2}\}}
    \def\qho{${\EuScript O}_{q,{\rm op}}$}
    \def\qhod{${\EuScript O}_{q,\ddc}$}
    \def\qhow{${\EuScript O}_{q,{\rm w}}$}

    \def\pdef{$\EuScript {PD}$}

    \theoremstyle{plain}

   \newtheorem{thm}{Theorem}
   \newtheorem{pro}[thm]{Proposition}
   \newtheorem{lem}[thm]{Lemma}
   \newtheorem{cor}[thm]{Corollary}

   \newtheorem*{samp}{Sample Theorem}

   \theoremstyle{definition}
   
   \newtheorem{exa}[thm]{{\it Example}}

   \theoremstyle{remark}
   \newtheorem{rem}[thm]{{\it Remark}}

\begin{document}

   \title[Operators of the $q$--oscillator  ]{Operators of the $q$--oscillator   }
   \author[F.H. Szafraniec]{Franciszek Hugon Szafraniec}
   \address{Instytut        Matematyki,         Uniwersytet
   Jagiello\'nski, ul. Reymonta 4, PL-30059 Krak\'ow}
   \email{fhszafra@im.uj.edu.pl}
\thanks{Supported at its final stage by the MNiSzW grant N201 026 32/1350}
   \begin{abstract}
    We scrutinize the possibility of extending the result of
\cite{ccr} to the case of $q$-deformed oscillator for $q$ real;
for this we exploit the whole range of the deformation parameter
as much as possible. We split the case into two depending on
whether a solution of the commutation relation is bounded or not.
Our {\it leitmotif}\/ is {\it subnormality}.

The deformation parameter $q$ is reshaped and this is what makes
our approach effective. The newly arrived parameter, the operator
$C$, has two remarkable properties: it separates in the
commutation relation the annihilation and creation operators from
the deformation as well as it $q$-commutes with those two. This is
why introducing the operator $C$ seems to be far-reaching.
   \end{abstract}
  \maketitle


$q$-deformations of the quantum harmonic oscillator (the
abbreviation the $q$-{\it oscillator} stands here for it) has been
arresting attention of many\,\footnote{\;$q$-deformations are
vastly disseminated in Mathematical Physics and we would like to
acknowledge here with pleasure \cite{kli} for bringing them closer
to Mathematics} resulting among other things in quantum groups.
Besides realizing the ever lasting temptation to generalize
matters, it brings forth new attractive findings. This paper
exhibits the {\em spatial} side of the story.

    The $q$-oscillator algebra, which is the {\it milieu} of our
    considerations, is that generated by three objects $a_+$,
    $a_-$ and $1$ (the latter being a unit in the algebra)
    satisfying the commutation relations
    \begin{equation}\label{w28.2}
a_-a_+-qa_+a_-=1;
    \end{equation}
    it goes back to the seventies with \cite{ari} as a specimen.
    The other versions which appear in the literature are
    equivalent to that and this is described completely in
    \cite{kli} where a list of further references can be found.

    Looking for $*$-representations of \eqref{w28.2} usually means
    assuming that $a_-=a_+^*$, with the asterisk denoting the
    Hilbert space adjoint. Thus what we start with is a {\it
    given} Hilbert space and the commutation relation
    \begin{equation} \tag{\qho}
    S\gw S-qSS\gw=I,
    \end{equation}
    in it. Of course, $q$ must be perforce \underline{real} then;
    this is what assume in the paper.

   An easy-going consequence is
    \begin{samp}
If $S$ is a weighted \label{b.1} shift with respect to the basis $
\{e_{n}\}_{n=0}^{\infty} $ and
    \begin{equation*}
S\gw Sf-qSS\gw f=f, \quad f\in\lin\{e_{n}\}_{n=0}^{\infty},
    \end{equation*}
    then $Se_n=\sqrt{1+q+\cdots+q^{n}}\,e_{n+1}$, $n\Ge 0$.
    \end{samp}
    `If $S$ is a weighted shift' -- this is usually tacitly
    assumed when dealing with the relation \naw{\qho}, like in
    \cite{che}. It is sometimes made a bit more explicit in
    stating that a vacuum vector (or a ground state, depending on
    denomination in Mathematical Physics an author belongs to) of
    $S$ exists. The point here (as it was in \cite{ccr} for $q=1$)
    is to discuss the case. It turns out that, like in \cite{ccr},
    {\it subnormality} plays an important role in the matter (and
    this, the case $q=1$ at least, is parallel to Rellich-Dixmier
    \cite{rel,dix} characterization of solutions to the CCR).
    Luckily, the above coincides with our belief that subnormality
    is the missing counterpart of complex variable in the
    quantization scheme.

    \section*{Preliminary essentials}
    \subsection*{A short guide to subnormality.} Recall
that a densely defined operator $A$ is said to be {\em
hyponormal}\/ if $\dz{A}\subset\dz{A\gw}$ and $\|A\gw
f\|\Le\|Af\|$, $f\in\dz A$. A hyponormal operator $N$ is said to
be {\it formally normal}\/ if $\|Nf\|=\|N\gw f\|$, $f\in\dz{N}$.
Specifying more, a formally normal operator $N$ is called {\em
normal}\/ if $\dz N=\dz{N\gw}$. Finally, a densely defined
operator $S$ is called {\em $($formally$)$ subnormal} if there is
a Hilbert space $\kkc$ containing $\hhc$ isometrically and a
(formally) normal operator $N$ in $\kkc$ such that $S\subset N$.

The following diagram relates these notions.
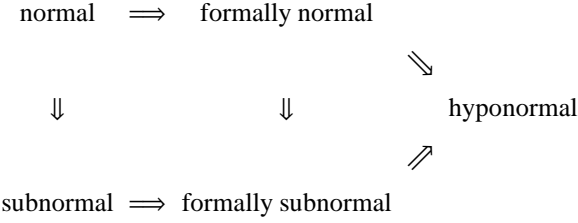
\setlength{\tabcolsep}{1mm}
   \renewcommand\arraystretch{1.5}
   \begin{equation*}
   \begin{tabular}{ccccc}
  normal&$\Longrightarrow$&formally normal&{}&{}
\\{}&{}&{}&{$\Searrow$}&{}
   \\ $\Downarrow$&{}&$\Downarrow$&{}&hyponormal \\{}&{}&{}&{$\Nearrow$}&{}
   \\ subnormal&$\Longrightarrow$&formally subnormal&{}&{}
   \end{tabular}
   \end{equation*}
   Though the definitions of formal normality and normality look
   much alike, with a little difference concerning the domains
   involved, the operators they define may behave in a totally
   incomparable manner. However, needless to say, these two
   notions do not differ at all in the case of bounded operators.

   If $A$ and $B$ are densely defined operators in $\hhc $ and
   $\kkc$ resp such that $\hhc\subset\kkc$ and $A\subset B$ then
   \begin{equation}\label{23.2.0}
   \dz A\subset \dz B\cap\hhc,\quad \dz{B\gw}\cap\hhc\subset
   P\dz{B\gw}\subset \dz{A\gw}
   \end{equation}
 where $P$ stands for the orthogonal projection of $\kkc$ onto
 $\hhc$; moreover,
   \begin{equation} \label{23.2.2}
   A\gw Px=PB\gw x, \quad x\in \dz {B\gw}.
   \end{equation}
   If $B$ closable, then so is $A$ and both $A\gw$ as well as
   $B\gw$ are densely defined. The extension $B$ of $A$ is said to
   be {\it tight}\/ if $\dz {\bar A}= \dz {\bar B}\cap\hhc$ and
   $*$-{\it tight}\/ if $\dz{B\gw}\cap\hhc=\dz{A\gw}$. If $\dz
   B\subset \dz {B\gw}$ (and this happens for formally normal
   operators as we already know), the two chains in \eqref{23.2.0}
   glue together as\,\footnote{\;Description of domains of
   weighted shifts and their adjoint can be found in \cite{ass}.}
   \begin{equation} \label{23.2.1}
   \dz A\subset \dz B\cap\hhc\subset\dz{B\gw}\cap\hhc\subset
   P\dz{B\gw}\subset \dz{A\gw} .
   \end{equation}
   As we have already said a densely defined operator having a
   normal extension is just {subnormal}. However, normal
   extensions may not be uniquely determined in unbounded case as
   their minimality becomes a rather fragile matter, see
   \cite{try3}; even though the inclusions \eqref{23.2.1} hold for
   any of them. Moreover, even if all of them turn into equalities
   none of the normal extensions may be minimal of cyclic type
   (this is what ensures uniqueness); this will become effective
   when we pass to the case of $q>1$. So far we have got an
   obvious fact.
   \begin{pro}\label{sroda}
   A subnormal operator $S$ has a normal extension which is both
   tight and $*$--tight if and only if
   \begin{equation} \label{21.1}
   \dz{\bar S}=\dz{S\gw}.
   \end{equation}
If this happens then any normal extension is both tight and
$*$--tight.
   \end{pro}
   Because equality \eqref{21.1} is undoubtedly decisive for a
   solution of the commutation relation of (any of) the
   oscillators to be a weighted shift, subnormality is properly
   settled into this context.

    \subsection*{$q$-notions.} For $x$ an integer and $q$ real,
$[x]_q \okr (1-q^{x})(1-q)^{-1}$ if $q\neq 1$ and $[x]_1\okr x $.
If $x$ is a non--negative integer, $[x]_q=1+q\cdots+q^{x-1}$ and
this is usually referred to as a {\it basic} or $q$--number. A
little step further, the $q$--factorial is like the conventional,
$[0]_q!\okr 1$ and $[n]_q!\okr[0]_q\cdots[n-1]_q[n]_q$ and so is
the $q$--binomial $\bbinom{m}{n}_q\okr
\frac{[m]_q!}{[m-n]_q![n]_q!}$. Thus, if $-1\Le q$ and
$x\in\mathbb N$ the basic number $[x]_q$ is non--negative.

For arbitrary complex numbers $a$ and $q$ one can always define
$(a;q)_k$ as follows
    \begin{equation*}
(a;q)_0\okr 1,\quad (a;q)_k\okr
(1-a)(1-aq)(1-aq^2)\cdots(1-aq^{k-1}),\quad k=1,2,3,\dots
    \end{equation*}
    Then for $n>0$ one has $[n]_q!=(q,q)_n(1-q)^{-n}$. Moreover,
    there are (at least) two possible definitions of
    $q$--exponential functions
   \begin{equation*}
e_q(z)\okr\sum_{k=0}^\infty\frac{1}{(q;q)_k}z^k,\quad
z\in\omega_q,
    \end{equation*}
    \begin{equation*}
E_q(z)\okr\sum_{k=0}^\infty\frac{q^{\binom{k}{2}}}{(q;q)_k}z^k,
\quad z\in\omega_{q^{-1}},\quad q\neq 0,
    \end{equation*}
    where
    \begin{equation*}
\omega_q\okr
    \begin{cases}
    \zb{z}{|z|<1} & \text{if $|q|<1$,}
\\
    \mathbb C& \text{otherwise}.
    \end{cases}  \smallskip
   \end{equation*}
   These two functions are related via
   \begin{equation*}
   e_q(z)=E_{q^{-1}}(-z),\quad z\in\omega_q,\quad q\neq 0.
    \end{equation*}

    \section*{The $q$ oscillator}
    \subsection*{Spatial interpretation of $\text{(\qho)}$.}

    The relation (\qho) has nothing but a symbolic meaning unless
    someone says something more about it; this is because some of
    the solutions may be unbounded. By reason of this we
    distinguish two, extreme in a sense, ways of looking at the
    relation \naw{\qho}:

   The first meaning of (\qho) is
    \begin{gather}\begin{split}
    \text{$S$ closable, $\ddc$ is dense in
 $\hhc$  and}\phantom{aaaaaaaaai}\\
\text{$\ddc\subset\dz{S\gw \bar S}\cap\dz{\bar SS\gw}$, $S\gw
Sf-qSS\gw f=f$, $f\in\ddc.$}\end{split}\tag{\qhod}
    \end{gather}
    The other is
    \begin{equation}
    \text{$\is{Sf}{Sg}-q\is{S\gw f}{S\gw g}=\is{f}{g}$, $f,g\in\dz
    S\cap\dz{S\gw}$}\tag{\qhow}
    \end{equation}
    and, because this is equivalent to
    \begin{equation*}
    \|Sf\|^2-q\|S\gw f\|^2=\|f\|^2,\quad f\in\dz S\cap\dz{S\gw}
    \end{equation*}
    it implies for $S$ to be closable, (\qhow) in turn is
    equivalent to
    \begin{equation*}
    \is{\bar Sf}{\bar Sg}-q\is{S\gw f}{S\gw g}=\is{f}{g},\quad
    f\in\dz {\bar S}\cap\dz{S\gw}.
    \end{equation*}

    The occurring interdependence, which follows, let us play
    variation on the theme of \naw{\qho}.
    \begin{enumerate}
    \item[\liczp 1] {\it  {\rm  (}\qhod{\rm  )}  with  $\ddc$
being a core of $S$ \wyn {\rm(}\qhow{\rm)} and $\dz{\bar
S}\subset\dz{S\gw}$.}
    \end{enumerate}
    \noindent Indeed, for $f\in\dz{\bar S}$ there is a sequence
    $(f_n)_n \subset\ddc$ such that $f_n\rightarrow f$ and
    $Sf_n\rightarrow \bar Sf$. Because $S\gw$ is closed we get
    from (\qhod) that $S\gw f_n\rightarrow S\gw f$ and
    consequently $f\in\dz{S\gw}$ as well as (\qhow).
    \begin{enumerate}
    \item[\liczp 2] {\it  {\rm  (}\qhod{\rm  )}  with  $\ddc$
being a core of $S\gw$ \wyn {\rm(}\qhow{\rm)} and $\dz{
S\gw}\subset\dz{\bar S}$.}
    \end{enumerate}
    \noindent This uses the same argument as that for \liczp 1.
    \begin{enumerate}
    \item[\liczp 3]  {\it  {\rm  (}\qhow{\rm  )}  \wyn  {\rm
(}\qhod{\rm )} with $\ddc=\dz{S\gw \bar S}\cap\dz{\bar SS\gw}$}.
    \end{enumerate}
    \noindent This is because $\dz{S\gw \bar S}\cap\dz{\bar
    SS\gw}\subset\dz{\bar S}\cap\dz{S\gw}$.
    \begin{enumerate}
    \item[\liczp 4] {\it {\rm (}\qhow{\rm )}  and  $\dz{\bar
S}\cap\dz{S\gw}$ a core of $S$ and $S\gw$ \wyn $\dz{S\gw\bar
S}=\dz{\bar SS\gw}$.}
    \end{enumerate}
    \noindent Take $f\in\dz{S\gw \bar S}$. This means $f\in
    \dz{\bar S}$ and $\bar Sf\in \dz{S\gw}$. Because of this,
    picking $(f_n)_n\in\dz{\bar S}\cap\dz{S\gw}$, we get from
    (\qhow) in limit
    \begin{equation}\label{x9}
    \is{S\gw\bar Sf}{g}-q\is{S\gw f}{S\gw g}=\is{f}{g}
    \end{equation} for $g\in\dz{\bar  S}\cap\dz{S\gw}$  and,
because $g\in\dz{\bar S}\cap\dz{S\gw}$ is a core of $S\gw$, we get
\eqref{x9} to hold for $g\in\dz{S\gw}$. Finally, $S\gw
f\in\dz{\bar S}$. The reverse inequality needs the same kind of
argument.

    The above results in
    \begin{enumerate}
    \item[\liczp 5] \label{piec}{\it  {\rm(}\qhow{\rm)}  and  $\dz{\bar
S}=\dz{S\gw}$ \wyn $\bar S$ satisfies {\rm(}\qhod{\rm)} on
$\ddc=\dz{S\gw\bar S}=\dz{\bar SS\gw}$.}
    \end{enumerate}
   \begin{rem} \label{tt1.4.4}
   Notice that when $q\neq -1$ and $S$ satisfying (\qhod) with
   $\ddc=\dz{S\gw\bar S}=\dz{\bar SS\gw}$ for $\ddc$ to be a core
   of $S\gw$ is necessary and sufficient $\ob {S\gw S}$ to be
   dense in $\hhc$.
   \end{rem}

   The following is a kind of general observation and settles
   hyponormality (or boundedness) in the context of (\qhod).
    \begin{pro}\label{t2}
    \naw{a} For $0\Le q<1$ and for $S$ satisfying \naw{\qhod},
    $S\res\ddc$ is hyponormal if and only if $S$ is bounded and
    $\|S\|\Le(1-q)^{-1/2}$. \naw{b} For $q<0$ and for $S$
    satisfying \naw{\qhod}, $S\gw\res\ddc$ is hyponormal if and
    only if $S$ is bounded and $\|S\|\Le(1-q)^{-1/2}$.
    \end{pro}
    \begin{proof}
    Write (\qhod) as
    $$
    (1-q)\|Sf\|^2=q(\|S\gw f\|^2-\|Sf\|^2)+\|f\|^2,\quad f\in\ddc.
    $$
    and look at this.
    \end{proof}

    \subsection*{The selfcommutator.}
    Assuming $\ddc\subset\dz{SS\gw}\cap\dz{S\gw S}$ we introduce
    the following operator
    \begin{equation}\label{1.1}
    C\okr I+(q-1)SS\gw,\quad \dz C\okr\ddc.
    \end{equation}
  This operator turns out to be an important invention in the
  matter. In particular there are two immediate consequences of
  this definition. The first says if $S$ satisfies (\qhod) with
  $\ddc$ invariant for both $S$ and $S\gw$ then $\ddc$ is
  invariant for $C$ as well and
    \begin{equation}\label{a2}
    CSf=qSCf,\quad qCS\gw f=S\gw Cf,\quad f\in\ddc.
    \end{equation}
 The other is that (\qhod) takes now the form
    \begin{equation}\label{a1}
    S\gw Sf-SS\gw f=Cf,\quad f\in\ddc,
    \end{equation}
  which means that $C$ is just the selfcommutator of $S$ on
  $\ddc$.

    We would like to know the instances when $C$ is a positive
    operator.
    \begin{pro}\label{t1}
    \naw{\rm a} For $q\Ge 1$, $C> 0$ always. \naw{\rm b} For
    $q<1$, $C\Ge 0$ if and only if $S$ is bounded and
    $\|S\|\Le(1-q)^{-1/2}$. \naw{\rm c} For $S$ satisfying
    \naw{\qhod}, $C\Ge 0$ if and only if $S$ is hyponormal.
    \end{pro}
    \begin{proof}
    While (a) is apparently trivial (b) comes out immediately from
    $$
    \is{Cf}{f}=\|f\|^2+(q-1)\|S\gw f\|^2,\quad f\in\ddc.
    $$
    For (c) write (using (\qhod)) with $f\in\ddc$
    $$
    \is{Cf}{f}=\|f\|^2+(q-1)\|S\gw f\|^2=\|f\|^2+q\|S\gw f\|^2
    -\|S\gw f\|^2=\| Sf\|^2-\|S\gw f\|^2.
    $$
    \end{proof}

    \begin{exa} \label{tw22.1}
    On the other hand, with any unitary $U$ the operator
    \begin{equation} \label{w22.1}
    S\okr (1-q)^{-1/2}U
    \end{equation}
     satisfies \naw{\qhod} if $q<1$. The operator $S$ is
     apparently bounded and normal. Consequently (the Spectral
     Theorem) it may have a bunch of nontrivial reducing subspaces
     (even not necessarily one dimensional) or may be irreducible
     and this observation ought to be dedicated to all those who
     start too fast generating algebras from formal commutation
     relations.
      \end{exa}
    \begin{pro} \label{tw28.1}
    For $q<1$ the only formally normal operators satisfying
    \naw{\qhod} are those of the form \eqref{w22.1}. For $q\Ge 1$
    there is no formally normal solution of \naw{\qhod}.
    \end{pro}
    \begin{proof} Straightforward.
   \end{proof}

    \begin{exa}\label{t4}
    An {\it ad hoc} illustration can be given as follows. Take a
    separable Hilbert space with a basis
    $(e_n)_{n=-\infty}^\infty$ and look for a bilateral (or rather
    {\it two-sided}) weighted shift $T$ defined as
    $Te_n=\tau_ne_{n+1}$, $n\in\mathbb Z$. Then, because $T\gw
    e_n=\bar\tau_{n-1}e_{n-1}$, $n\in\mathbb Z$, for any
    $\alpha\in\mathbb C$ and $N\in\mathbb Z$ we get
    $|\tau_n|^2=\alpha q^{n+N}+(1-q^{n+N})(1-q)^{-1}=\alpha
    q^{n+N}+[n+N]_q$ for all $n$ if $q\neq 1$ and
    $|\tau_n|^2=\alpha+n$ if $q=1$; this is for all $n\in\mathbb
    Z$. The only possibility for the right hand sides to be
    non--negative (and in fact positive)\,footnote{\;We avoid
    weights which are not non--negative, for instance complex, as
    they lead to a unitary equivalent version only.} is $\alpha\Ge
    (1-q)^{-1}$ for $0\Le q<1$ and $\alpha=(1-q)^{-1}$ for $q<0$;
    the latter corresponds to Example \ref{w22.1}. Thus {\it the
    only} bilateral weighted shifts satisfying (\qhod), with
    $\ddc=\lin\zb{e_n}{n\in\zzz}$, are those $Te_n=\tau_ne_{n+1}$,
    $n\in\mathbb Z$ which have the weights
   \begin{equation*}
   \tau_n\okr
   \begin{cases}
   \sqrt{(1-q)^{-1}},& \phantom{0\Le}\;q\Le0\\
    \sqrt{\alpha q^{n+N}+[n+N]_q},\text{\quad
$\alpha>(1-q)^{-1}$,\; $N\in\mathbb Z$,\quad} &0\Le q<1\\
   \phantom{l}\text{none},& 1\Le q
   \end{cases}
   \end{equation*}
   However, $T$ {\it violates hyponormality} (pick up $f=e_0$ as a
   sample) if $0<q<1$. Also $C$ defined by \eqref{1.1} is {\it
   neither positive nor negative} ($\is{Ce_0}{e_0}=a>0$ while
   $\is{Ce_{-1}}{e_{-1}}<0$). Let us mention that $T$ is
   $q^{-1}$--hyponormal in the sense of \cite{ota}. Anyway, $T$ is
   apparently {\it unbounded} if $q>0$. The case of $q\Le 0$ is
   precisely that of Example \ref{w22.1}.
    \end{exa}

   \begin{exa} \label{t10.3.1t}
   Repeating the way of reasoning of Example \ref{t4} we get that
   {\it the only} unilateral weighted shifts satisfying (\qhod)
   are those $T$, defined as $Te_n=\tau_ne_{n+1}$ for $n\in\mathbb
   N$, which have the weights
   \begin{equation*}
   \tau_n=\sqrt{\,[n+1]_q}\;,\quad -1\Le q.
   \end{equation*}
   This is so because the virtual, in this case, $`\tau_{-1}{}'$
   is $0$ ($T\gw e_0=0$). If $-1\Le q<0$ they are bounded and {\it
   not hyponormal}, if $0\Le q<1$ they are again bounded and {\it
   hyponormal} and if $1\Le q$ the are unbounded and {\it
   hyponormal}; the two latter are even subnormal (cf. Theorem
   \ref{t5.1} and \ref{t27.2.1t} resp.).
   \end{exa}

   \begin{rem} \label{tg1}
   According to Lemma 2.3 of \cite{kii} for $0<q<1$ the only cases
   which may happen are the orthogonal sums of the operators
   considered in Examples \ref{t4}, \ref{t10.3.1t} and given by
   formula \eqref{w22.1}. For $q>1$, due to the same Lemma, the
   orthogonal sum of that from Example \ref{t10.3.1t} can be taken
   into account.
   \end{rem}
    \subsection*{An auxiliary lemma of \cite{sch1}.}
    We state here a result, \cite{sch1} Lemma 2.4, which
    authorizes the examples above. We adapt the notation of
    \cite{sch1} to ours as well as improve a bit the syntax of the
    conclusion therein.
   \begin{lem} \label{tsch1}
   Let $0<p<1$ and $\varepsilon\in\{-1,+1\}$. Assume $T$ is a
   closed densely defined operator in $\hhc$. Then
   \begin{equation} \label{schi}
   T\gw Tf-p^2TT\gw f=\varepsilon(1-p^2)f,\quad f\in\dz{T\gw
   T}=\dz{TT\gw}
   \end{equation}
 if and only if $T$ is unitarily equivalent to an orthogonal sum
 of operators of the following type:

   \noindent $\boldsymbol \cdot\,$ in the case of $\varepsilon=1$
   \begin{enumerate}
   \item[{\rm (I)}] $T_{\rm I}\colon f_n\to(1-p^{2(n+1)})^{1/2}f_{n+1}$ in $\hhc=\bigoplus\nolimits_{n=0}
^{+\infty}\hhc_n$ with each $\hhc_n\okr\hhc_0$;
   \item[{\rm (II)}] $T_{\rm II}\colon f_n\to(1+q^{2(n+1)}A^2)^{1/2}f_{n+1}$ in $\hhc=\bigoplus\nolimits_{n=-\infty}
^{+\infty}\hhc_n$ with each $\hhc_n\okr\hhc_0$ and $A$ being a
selfadjoint operator in $\hhc_0$ with ${\rm sp}(A)\subset[p,1]$
and either $p$ or $1$ not being an eigenvalue of $A$;
   \item[{\rm (III)}] $T_{\rm III}$ a unitary operator;
   \end{enumerate}
   \noindent $\boldsymbol \cdot\,$ in the case of $\varepsilon=-1$
   \begin{enumerate}
   \item[{\rm (IV)}]  $T_{\rm IV}\colon f_n\to(p^{2n}-1)^{1/2}f_{n-1}$ in $\hhc=\bigoplus\nolimits_{n=0}
^{+\infty}\hhc_n$ with each $\hhc_n\okr\hhc_0$ and always
$f_{-1}\okr0$.
   \end{enumerate}
   \end{lem}
   A couple of remarks seem to be absolutely imperative.
   \begin{rem} \label{trsch1}
   The conclusion of Lemma \ref{tsch1} is a bit too condensed. Let
   us provide with some hints to reading it. First of all the way
   of understanding the meaning of $f_n$'s appearing in (I), (II)
   and (IV) should be as follows: take $f\in\hhc_0$ and define
   $f_n$ as a (one sided or two sides, depending on circumstances)
   sequence having all the coordinates zero except that of number
   $n$ which is equal to $f$. Then, with a definition
$$\EuScript D(\mathcal
E)\okr\lin\zb{f_n}{f\in\mathcal E\subset\hhc_0,\;n\in\mathbb Z
\text{ or }n\in\mathbb N\text{ depending on the case}},$$ one has
to guess that $\dz{T_{\rm I}}=\dz{T_{\rm IV}}=\EuScript D(\hhc_0)$
and $\dz{{T_{\rm II}}}=\EuScript D(\dz A)$. Passing to closures in
(I), (II) and (IV) we check that $\overline {T}_{\rm I}$ as well
as $\overline{ T}_{\rm IV}$ are everywhere defined bounded
operators (use $0<p<1$) while $\overline{ T}_{\rm II}$ is always
unbounded (though satisfying $\dz{{T}_{\rm
II}^{\,*}\overline{T}_{\rm II}}=\dz{\overline{T}_{\rm II}{T}_{\rm
II}^{\,*}}$\,\footnote{\;In this matter we have implications
\liczp 4 and \liczp 5 on p. \pageref{piec}.}).
   \end{rem}

   \begin{rem} \label{trsch2}
    To relate \eqref{schi} to (\qhod) set $\varepsilon=1$,
    $p=\sqrt q$ and $T=\sqrt{1-p^2}S$ when $0<q<1$, and
    $\varepsilon=-1$, $p^{-1}=\sqrt q$ and
    $T=p^{-1}\sqrt{p^2-1}S\gw$ when $q>1$.
   \end{rem}

    \subsection*{Positive definiteness from (\qhod).}
   The following formalism will be needed.
    \begin{pro} \label{tx}
     If $S$ satisfies {\rm(}\qhod{\rm)} with $\ddc$ being
     invariant for both $S$ and $S\gw$, then
    \begin{equation}\label{x5}
    S\gw{}^iS^jf=\sum_{k= 0}^\infty
    [k]_q!\bbinom{i}{k}_q\bbinom{j}{k}_qS^{j-k} C^k
    S\gw{}^{(i-k)}f,\quad f\in\ddc,\; i,j=0,1,\dots,
    \end{equation}
    If, moreover, $C\Ge 0$ then
    \begin{equation}\label{x8}
    \sum_{i,j= 0}^ p\is{S^if_j}{S^jf_i}=\sum_{k =0}^\infty
    [k]_q!\,\left\|\sum_{i= 0}^
    p\bbinom{i}{k}_qC^{k/2}S\gw{}^{(i-k)}f_i\right\|^2,\quad
    f_0,\dots f_p\in\ddc.
    \end{equation}
    All this under convention $S^l=(S\gw)^l=0$ for $l<0$ and
    $\bbinom{i}{j}_q=0$ for $j>i$.
    \end{pro}
    \begin{proof}
    Formula \eqref{x5} is in \cite[formula (35)]{cig}. Formula
    \eqref{x8} is an immediate consequence of \eqref{x5}.
    \end{proof}

    As a direct consequence of Fact A and \eqref{x8} we get
    \begin{cor} \label{tw22.3}
    Suppose $S$ satisfies {\rm(}\qhod{\rm)} with $\ddc$ being
    invariant for $S$ and $S\gw$ as well as $\ddc$ is a core of
    $S$. If $C\Ge 0$, then
    \begin{equation}
    \sum_{i,j= 0}^ p\is{S^if_j}{S^jf_i}\Ge 0,\quad f_0,\dots
    f_p\in\ddc. \tag{\pdef}
    \end{equation}
    \end{cor}

   \subsection*{A useful Lemma.}
   \begin{lem} \label{tt1.1.g}  Let $q>0$. Consider following conditions:
   \begin{enumerate}
   \item[{\rm (a)}]  $S$  satisfies  {\rm  (}\qhow{\rm   )}   and
$\dz{\bar S}=\dz{S\gw}$;
   \item[{\rm (b)}]  $\jd{S\gw}\neq\nul$ and   for $n=0,1,\dots$
   \begin{equation}\label{g17}
   \text{$f\in\jd{S \gw}$ \wyn $\bar S ^n f\in\dz{\bar S }$,\;
   $\bar S^ {(n-1)}f\in\dz{S \gw}$\, \& \;$S \gw\bar S
   ^{n-1}f=(n-1)\bar S^{n-2}f$};
   \end{equation}
   \item[{\rm (c)}] there is $f\neq0$ such that $\bar S^nf\in\dz{\bar S}$,
$n=0,1,\dots$ and $\bar S^mf\perp \bar S^n$ for $m\neq n$.
   \end{enumerate}
   Then {\rm (a)} \wyn {\rm (b)} \wyn {\rm (c)}.
   \end{lem}
   \begin{proof}
   (a) \wyn (b). The polar decomposition for $S\gw$ is
   $S\gw=V|S\gw|$ where $V$ is a partial isometry with the initial
   space $\ob{|S\gw|}$ and the final space $\ob{SS\gw}$. Suppose
   $\jd{S\gw}=\nul$. Then, because $\jd
   V=\ob{|S\gw|}^\perp=\jd{|S\gw|}=\jd{\bar SS\gw}=\jd{S\gw}$, $V$
   is unitary. Since $\bar S=|S\gw|V\gw$, from \liczp 5 we get
   $V|S\gw|^2V\gw=q| S\gw|^2+I$. Consequently, for the spectra we
   have ${\rm sp}(|S\gw|)\subset q\,{\rm sp}({|S\gw|})+1\subset
   [0,+\infty)$ which is an absurd. Thus $\jd{S\gw}\neq\nul$.

   We show \eqref{g17} by induction. Of course, $\jd
   {S\gw}\subset\dz {\bar S }=\dz{S \gw}$, which establishes
   \eqref{g17} for $n=0$. Suppose $\jd {S \gw}\subset\dz{\bar S
   ^n}$ and $S \gw\bar S ^{n-1}f=(n-1)\bar S ^{n-2}f$. Then, for
   $g\in\dz {\bar S}=\dz{S \gw}$,

   \begin{equation}\label{g18}
   \is{S \gw\bar S ^{n-1}f}{S \gw g}=(n-1)\is{\bar S ^{n-2}f}{\bar
   S\gw g}.
   \end{equation}
   Because already $\bar S ^{(n-2)}f\in\dz{\bar S
   }=\dz{S\gw{}\gw}$, we have

   \begin{equation}\label{g19}
   |\is{S \gw \bar S ^{n-1}f}{S \gw g}|\le C\|g\|.
   \end{equation}
   Because $\bar S ^{(n-1)}\in\dz{\bar S }=\dz{S \gw}$, we can use
   (\qhow) so as to get
   \begin{equation*}
   \is{\bar S ^nf}{\bar S g}=\is{\bar S \bar S ^{(n-1)}f}{\bar S
   g}= \is{S \gw \bar S ^{(n-1)}}{S \gw}+\is{\bar S ^{(n-1)f}}{g}.
   \end{equation*}
   This, by \eqref{g19}, implies $\bar S ^n f\in\dz{S
   \gw}=\dz{\bar S }$ and, consequently, by \eqref{g18}, gives us
   $S \gw\bar S ^{n}f=n\bar S ^{n-1}f$ which completes the
   induction argument. Now a straightforward application of
   \eqref{g17} gives $\bar S^n(\jd{S\gw)} \subset \dz{\bar S}\cap
   \dz{S\gw}$ for $n=0,1,\dots$.

   (b) \wyn (c). Take any $f\in\jd{S\gw}$ and using \eqref{g17}
   and \eqref{x5} write
   \begin{align*}
   \is{S^m f}{S^nf }=\is{S^n{}\gw S^mf}{f}=\sum_{k= 0}^
   {\min\{m,n\}}[k]_q!\,\bbinom m k_q\, \bbinom n k_q\,
   \is{S^{(n-k)}C^kS\gw{}^{(m-k)}f }{f}=0,\quad m>n.
   \end{align*}
   \end{proof}

   \subsection*{A matrix formation.}  Suppose $q>0$ and $S$ is a weighted
shift with respect to $(e_k)_{k=0}^\infty$ with the weights
$(\!\sqrt{\,[k+1]_q}\,)_{k=0}^\infty$. With
   \begin{equation} \label{t1.21.2.7}
   S_0\okr S,\quad S_n \okr q^{n/2}S,\quad
   D_n\okr\sqrt{[n]_q}\,\diag(q^{k/2})_{k=0}^\infty,\quad
   n=1,2\dots
   \end{equation}
   the matrix
   \smallskip
   \begin{equation} \label{2.21.2.7}
    \left(\begin{array}{ccccc}
  S_{0} & D_{1} & 0 & 0 &{}\\
  0 & S_{1} &D _{2} & 0&\ddots \\
  0 & 0 & S_2 & D_3&\ddots \\
  {} &\ddots& \ddots & \ddots & \ddots
   \end{array}\right)
   \end{equation}
    \smallskip
    defines an operator ${\boldsymbol N}$ in
    $\bigoplus_{n=0}^\infty\mathcal H_n$, $\mathcal H_n=\hhc$,
    with domain composed of all those $\bigoplus_{n=0}f_n$ for
    which $f_n=0$ but a finite number of $n$'s. This matrix, for
    the familiar creation operator was set out in \cite{born}.

   First we need to determine $\dz{{\boldsymbol N}\gw}$ and relate
   it to $\dz{\boldsymbol N}$. If $0<q<1$ then each $D_n$ is
   bounded. In that case Remark 9 in \cite{try4} gives us
   \begin{equation} \label{g1}
   \dz{{\boldsymbol N}\gw}=\bigoplus_{n=0}^\infty\dz{S_n\gw}.
   \end{equation}
If $q>1$ then each $S_nD_n^{-1}$ is bounded. According to
Proposition 4.5 in \cite{man} and Corollary 8 in\cite{try4} we can
deduce \eqref{g1} as well. In either case, what we get is the
adjoint of $\boldsymbol N$ can be taken as a matrix of adjoints
(which is rather an exceptional case). Because the same argument
concerning the adjoint of a matrix operator applies now to
${\boldsymbol N}\gw$ we can assert that the closure operation for
the operator $\boldsymbol N$ goes entrywise as well. Now, due to
the fact that the apparent norm equality for ${\boldsymbol N}$ and
${\boldsymbol N}\gw$ holds on $\dz{{\boldsymbol N}}$, we get
essential normality of ${\boldsymbol N}$. Consequently,
   \begin{equation} \label{g2}
  \text{$S$ is subnormal and $\bar N$ is its tight and $*$--tight
  normal extension.}
   \end{equation}

    \section*{Subnormality in the $q$-oscillator}
    \subsection*{The case of $S$ bounded.}

    The next says a little bit more about boundedness of solutions
    of (\qhod).

    \begin{pro} \label{tw22.4}
    Suppose $S$ is bounded and satisfies \naw{\qhod}. \naw{a} If\/
    $q<0$ then $\|S\|\Ge (1-q)^{-1/2}$. \naw{b} If\/ $0\Le q<1$
    then $\|S\|\Le (1-q)^{-1/2}$. \naw{c} If\/ $q\Ge 1$ then no
    such an $S$ exists.
    \end{pro}
    \begin{proof}
For (a) look at $\|Sf\|^2= \|f\|^2+q\|S\gw f\|^2\Ge
\|f\|^2+q\|S\|^2\|f\|^2$, for (b) do at $\|Sf\|^2= \|f\|^2+q\|S\gw
f\|^2\Le \|f\|^2+q\|S\|^2\|f\|^2$. For (c) write $\|Sf\|^2=
\|f\|^2+q\|S\gw f\|^2\Ge q\|S\|^2\|f\|^2$ which gives $1\Ge q$.
The case of $q=1$ is excluded by the well known result of Winter.
    \end{proof}
    \subsubsection*{The case of\/ $q<0$.}

Here we get at once
    \begin{cor} \label{tw29.1}
    For $q<0$ the only bounded operator $S$ with norm $\|S\|=
    (1-q)^{-1/2}$ satisfying \naw{\qhod} is that given by
    \eqref{w22.1}.
    \end{cor}

    \begin{proof}
By Proposition \ref{tw22.4} (a) and Proposition \ref{t2} (b)
$S\gw\res\ddc$ is hyponormal. On the other hand, by Proposition
\ref{t1} (b) and (c) $S\res\ddc$ is hyponormal too. Proposition
\ref{tw28.1} makes the conclusion.
    \end{proof}
    {Pauli} matrices, which are neither hyponormal nor
    cohyponormal\,\footnote{\;An operator $A$ is said to be {\it
    cohyponormal} if $A\gw$ is hyponormal; for unbounded $A$ this
    may not be the same as $A\gw\res{\dz A}$ to be hyponormal.},
    provide an example of operators satisfying (${\EuScript
    O}_{-1,{\rm op}}$) with norm $1>2^{-1/2}=(1-q)^{-1/2}$. Are
    there bounded operators satisfying (\qho) with norm not to be
    equal $(1-q)^{-1/2}$ for arbitrary $q<0$, different from $-1$
    say?

    \subsubsection*{The case of\/ $0\Le q<1$.}

    We list two results which hold in this case
    \begin{pro}\label{t5}
    Suppose $S$ satisfies \naw{\qhod} with $\ddc$ dense in $\hhc$.
    If\/ $0\Le q<1$, then the following facts are equivalent
    \begin{enumerate}
    \item[{\rm(i)}] $S$ is bounded and $\|S\|\Le (1-q)^{-1/2}$;
   \item[{\rm(ii)}] $S$ is bounded; \item[{\rm(iii)}] $S$ is subnormal;
   \item[{\rm(iv)}] $S$ is hyponormal.
    \end{enumerate}
     \end{pro}

    \begin{proof}
Because of conclusion (a) of Proposition \ref{t1} the only
remaining implication to argue for is (ii) $\Rightarrow$ (iii).
But, in virtue of \eqref{x8}, this comes out from the Halmos-Bram
characterization \cite{bra} of subnormality of bounded operators.
    \end{proof}
    \begin{thm}\label{t5.1}
    If\/ $0\Le q<1$, then the following facts are equivalent
    \begin{enumerate}
    \item[{\rm(i)}] there is an orthonormal basis $(e_n)_{n=0}^\infty$ in $\hhc$ such that $Se_n=\sqrt{[n+1]_q} \,e_{n+1}$,
$n=0,1,\dots$;
    \item[{\rm(ii)}] $S$ is irreducible\,\footnote{\label{t} Let us recall relevant
     definitions: a subspace $\ddc\subset\dz A$ is {\it invariant}
     for $A$ if $A\ddc\subset\ddc$; $A\res\ddc$ stands for the
     restriction of $A$ to $\ddc$. On the other hand, a {\it
     closed} subspace $\mathcal L$ is {\it invariant} for $A$ if
     $A(\mathcal L\cap\dz A) \subset\dz A$; then the restriction
     $A\rres{\mathcal L}\okr A\res{\mathcal L\cap\dz A}$. A step
     further, a closed subspace ${\mathcal L}$ {\it reduces} an
     operator $A$ if both $\mathcal L$ and $\mathcal L^\perp$ are
     invariant for $A$ as well as $P\dz A\subset\dz A$, where $P$
     is the orthogonal projection of $\tilde\hhc$ onto $\mathcal
     L$; all this is the same as to require $P A\subset AP$. Then
     the restriction $A\rres{\mathcal L}$ is called a {\it part}
     of $A$ in $\mathcal L$. $A$ is {\it irreducible } if it has
     no nontrivial reducing subspace. Comparing to the more
     familiar case of bounded operators some nuances become
     requisite here. Therefore, if $\mathcal L$ reduces $A$, then
     $\overline{(A\rres{\mathcal L})}=\bar A\rres{\mathcal L}$ and
     $(A\rres{\mathcal L})\gw=A\gw\rres {\mathcal L}$}, satisfies
     \naw{\qhod} with some $\ddc$ dense in $\hhc$, is bounded and
     $\|S\|= (1-q)^{-1/2}$;
    \item[{\rm(iii)}] $S$ is irreducible, satisfies
     \naw{\qhod} with some $\ddc$ dense in $\hhc$, is bounded and
     $\|S\|\Le (1-q)^{-1/2}$;
    \item[{\rm(iv)}] $S$ is irreducible, satisfies \naw{\qhod} with some
$\ddc$ dense in $\hhc$ and is bounded;
   \item[{\rm(v)}] $S$ is
    irreducible, satisfies with some $\ddc$ dense in $\hhc$
    \naw{\qhod} and is subnormal;
   \item[{\rm(vi)}]
    $S$ is irreducible, satisfies \naw{\qhod} with some $\ddc$
    dense in $\hhc$ and is hyponormal.
    \end{enumerate}
    \end{thm}
    \begin{proof}
   Proposition \ref{t5} establishes the equivalence of (ii) up to
   (vi).

Because $\sup\{\sqrt{[n+1]_q};\;n\Ge 0\}=(1-q)^{-1}$ and for $S$
as being a weighted shift $\|S\|= \sup\{\sqrt{[n+1]_q};\;n\Ge
0\}$, we get (i) $\Rightarrow$ (ii).

    Assume (iv). Because $\dz{\bar S}=\dz{S\gw}$, condition (c) of
    Lemma \ref{tt1.1.g} let calculate the weights of $\bar S$
    while starting with $e_0\in\jd{N\gw}$. Because $S$ is
    irreducible the sequence $(e_n)_{n=0}^\infty$ is complete.
    This establishes (i).
    \end{proof}

    \begin{rem} \label{tw22.5}
    From Theorem \ref{t5.1} and Example \ref{tw22.1} we get that
    there are two, of different nature, solutions of (\qhod). Is
    there any other at all?
    \end{rem}

    \subsubsection*{The case of\/ $q>1$.}   No bounded solution exits at all, cf.
    Proposition \ref{tw22.4} part (c).

   Let us memorize what is known already in the bounded case by
   the following tableau.
   \setlength{\tabcolsep}{2mm}
   \renewcommand\arraystretch{1.5}
   \begin{center}
   \begin{tabular}
   {l|l|c|c|c} {}&{}&$q<0$&$0\Le q<1$&$1\Le q$\\\cline{1-5}normal
   &{general}& \shortstack{{}\\{}\\SOME\\ \scriptsize{Exa.
   \ref{w22.1}}}&{} \shortstack{{}\\{}\\SOME\\ \scriptsize{Exa.
   \ref{w22.1}}}&{}\\\cline{1-4}{}&unilat.
   shift&{}&\shortstack{{}\\{}\\SOME\\ \scriptsize{Th.
   \ref{t5.1}}} &{}\\\cline{2-4} subnormal&bilat. shift&
   \shortstack{{}\\{}\\NONE\\ \scriptsize{Exa. \ref{t4}}}&
   \shortstack{{}\\{}\\NONE\\ \scriptsize{Exa.
   \ref{t4}}}&{}\\\cline{2-4} {}&others& \shortstack{{}\\{}\\SOME\\
   \scriptsize{Exa. \ref{tw22.1}}}& \shortstack{{}\\{}\\SOME\\
   \scriptsize{Exa. \ref{tw22.1}}}&\shortstack{{}\\{}\\NONE \\
   \scriptsize{Prop. \ref{tw22.4}(a)}}\\\cline{1-4}{}&unilat.
   shifts&{}& \shortstack{{}\\{}\\SOME\\ \scriptsize{Th.
   \ref{t5.1}}}&{}\\\cline{2-4} hyponormal &bilat. shift&
   \shortstack{{}\\{}\\NONE\\ \scriptsize{Exa. \ref{t4}}}&
   \shortstack{{}\\{}\\NONE\\ \scriptsize{Exa.
   \ref{t4}}}&{}\\\cline{2-4} {}&other& \shortstack{{}\\{}\\SOME\\
   \scriptsize{Exa. \ref{tw22.1}}}& \shortstack{{}\\{}\\SOME\\
   \scriptsize{Exa. \ref{tw22.1}}}&{}
   \end{tabular}
   \end{center}

    \subsection*{The case of $S$ unbounded.}
\subsubsection*{The case of $q<0$.}  There is no hope to look for
subnormal solutions of (\qho) among weighted shifts, neither one-
nor two-sided.

The only one-sided weighted shifts satisfying (\qho) are for
$-1<q<0$ and they are given as in (i) of Theorem \ref{t5.1}. They
are apparently not hyponormal (their weights are not increasing).

   The only two-sided weighted shifts which satisfy (\qho) are
   those of Example \ref{t4}. They are normal bilateral weighted
   shifts. So if there are subnormal operators satisfying (\qho)
   they must not be weighted shifts or bounded operators of norm
   less or equal $(1-q)^{-1/2}$, cf. Corollary \ref{tw29.1}.

   \subsubsection*{The case of\/ $0\Le q<1$.}
Lemma \ref{tsch1} does not leave any hope subnormal solutions
different than those in Theorem \ref{t5.1} but they must
necessarily be bounded.
\subsubsection*{The case of\/ $q\Ge 1$.}    This is the right case
for unbounded solutions to exist.
\begin{thm}\label{t27.2.1t}
   For a densely defined closable operator $S$ in a complex
   Hilbert space $\hhc$ consider the following conditions
      \begin{enumerate}
      \item[{\rm(i)}] $\hhc$ is separable and there is an orthonormal
   basis in it of the form $\{e_{n}\}_{n=0}^\infty$ contained in
   $\dz {\bar S}$ and such that
   \begin{equation}\label{27 2.1}
\bar Se_{n}=\sqrt{\,[n+1]_q}\,e_{n+1},\quad n=0,1,\dots;
   \end{equation}
      \item[{\rm(ii)}] $S$ is irreducible, satisfies {\rm (}\qhod{\rm )} with some
   $\ddc$ being invariant for $S$ and $S\gw$ and being a core of
   $S$, and $S$ is a {\bf subnormal} operator having a tight and
   $*$-tight normal extension;
      \item[{\rm(iii)}] $S$ is irreducible, satisfies {\rm (}\qhod{\rm )} with some
   $\ddc$ being a core of both $S$ and $S\gw$;
   \item[{\rm(iv)}] $S$   is irreducible,
   satisfies {\rm (}\qhow{\rm )} and $\dz{\bar S}=\dz{S\gw}$;
      \item[{\rm(v)}] $S$ is irreducible, satisfies {\rm (}\qhow{\rm )} with $\dz{\bar
   S}\cap\dz{S\gw}$ being dense in $\hhc$, $\jd{S\gw}\neq\nul$ and
   $\bar S^n(\jd{S\gw)} \subset \dz{\bar S}\cap \dz{S\gw}$ for
   $n=0,1,\dots$.
   \end{enumerate}
   Then {\rm (i)} \wyn {\rm (ii)} \wyn {\rm (iii)} \wyn {\rm (iv)}
   \wyn {\rm (v)} \wyn {\rm(i)}.
   \end{thm}
   \begin{proof}
   The implication {\rm (i)} \wyn {\rm (ii)} comes out from
   \eqref{g2}. Proposition \ref{sroda} leads us from (ii) to
   (iii), from there using Lemma \ref{tt1.1.g} comes it up to (v).
   Now, like in the proof of Theorem \ref{t5.1}, calculating the
   weights rounds up the chain of implications.
  \end{proof}

Now we visualize this section findings in the following tableau.
   \setlength{\tabcolsep}{2mm}
\renewcommand\arraystretch{1.5}
   \begin{center}
   \begin{tabular}
   {l|l|c|c|c} {}&{}&$q<0$&$0\Le q<1$&$1\Le q$\\\cline{1-5}normal
   &general& {}&{} {}&{\shortstack{{}\\{}\\NONE\\
   \scriptsize{Prop.
   \ref{tw28.1}}}}\\\cline{1-2}\cline{5-5}{}&unilat. shift&{}&{}
&\shortstack{{}\\{}\\{\phantom{E}}\\
\scriptsize{{\phantom {E}}}}\\\cline{2-2}\cline{5-5}
subnormal&bilat. shift& {}& {}&\shortstack{{}\\{}\\NONE\\
\scriptsize{Exa. \ref{t4}}}\\\cline{2-2}\cline{5-5} {}&others&
\shortstack{{}\\{}\\NONE\\ \scriptsize{Prop. \ref{t2}(b)}}&
\shortstack{{}\\{}\\NONE\\ \scriptsize{Prop.
\ref{t2}(a)}}&{}\\\cline{1-2}\cline{5-5}{}&unilat. shifts&{}&
{}&\shortstack{{}\\{}\\{{\phantom {E}}}\\ \scriptsize{{{\phantom
{E}}}}}
\\\cline{2-2}\cline{5-5} hyponormal &bilat. shift&
{}& {}&\shortstack{{}\\{}\\NONE\\ \scriptsize{Prop.
\ref{t2}(b)}}\\\cline{2-2}\cline{5-5} { }&others& {}&
{}&\shortstack{{}\\{}\\MAY\\ \scriptsize{Prop. \ref{t1}(a)\&(b)} }
   \end{tabular}
   \end{center}

    \section*{The $q$ oscillator: models in RKHS}
   \subsection*{A general look at.}    A
reproducing kernel Hilbert space $\hhc$ and its kernel $K$ which
suits our considerations is of the form
    \begin{equation}
\label{1.6.5.7} K(z,w)\okr \sum_{n=0}^{+\infty}
c_nz^n\overline{w}^n, \quad z,w\in D, \quad D=\mathbb C \text{ or
} D=\zb{z}{|z|<R\Le 1}.
   \end{equation}
   Notice $(\sqrt {c_n}Z^n)_{n=0}^{+\infty}$ is an orthonormal
   basis of $\hhc$.

   The following fact comes out, as a byproduct, from some general
   results on subnormality in \cite{ss2}; we give here an {\it ad
   hoc} argument. Let us make a shorthand notation
   \begin{equation} \label{11.14.5.7}
   \text{$\hhc\subset\Ld{\mathbb C,\mu}$ isometrically.}
   \end{equation}

   \begin{pro}\label{t1.12.5.7}
   There is a measure $\mu$ such that \eqref{11.14.5.7} holds if
   and only if there is a Stieltjes moment sequence
   $(a_n)_{n=0}^{+\infty}$ such that
   \begin{equation} \label{3.12.5.7}
    a_{2n}=c_{n}^{-1},\quad n=0,1,\dots
   \end{equation}
   If this happens than a measure $\mu$ can be chosen to be
   rotationally invariant\,\footnote{\;Or {\it radial} as some
   authors say.}, that is such that $\mu(\E^{\I
   t}\!\sigma)=\mu(\sigma)$ for all $t$'s and $\sigma$'s.
   \end{pro}
  \begin{proof}
   Suppose \eqref{11.14.5.7} to hold. Because $(\sqrt
   {c_n}Z^n)_{n=0}^{+\infty}$ is an orthonormal sequence in
   $\Ld{\mathbb C,\mu}$, we have
   \begin{equation*}
   c_n\od=\int_{\mathbb C}|z|^{2n}\mu(\D z),\quad n=0,1,\dots
   \end{equation*}
   Let $m_\mu$ be the measure on $[0,+\infty)$ transported from
   $\mu$ via the mapping $\mathbb C\ni z\to |z|\in[0,+\infty)$.
   Then
   \begin{equation}\label{2.12.5.7}
   a_n\okr\int_0^{+\infty}r^nm_\mu(\D r)=\int_{\mathbb
   C}|z|^{n}\mu(\D z),\quad n=0,1,\dots
   \end{equation}
    satisfies \eqref{3.12.5.7} as well as the sequence
    $(a_n)_{n=0}^{+\infty}$ is a Stieltjes moment sequence.

   If $(a_n)_{n=0}^{+\infty}$ is any Stieltjes moment sequence
   with a representing measure $m$ and satisfying \eqref{3.12.5.7}
   then the rotationally invariant measure
   \begin{equation}\label{1.13.5.7}
\mu(\sigma)\okr(2\pi)^{-1}\int_0^{2\pi}\int_0^{+\infty}\chi_\sigma(r\E^{\I
t})m(\D r)\D t,\quad \sigma \text{ Borel subset of } \mathbb C
   \end{equation}
   makes the imbedding \eqref{11.14.5.7} happen.
   \end{proof}

   \begin{thm} \label{t3.12.5.7}
  Under the circumstances of\/ {\rm Proposition \ref{t1.12.5.7}}
  there exists a not rotationally invariant measure $\mu$ such
  that \eqref{11.14.5.7} holds if and only if there is a sequence
  $(a_n)_{n=0}^{+\infty}$ satisfying \eqref{3.12.5.7} which is not
  Stieltjes determinate.
   \end{thm}
   \begin{proof}
   Suppose \eqref{11.14.5.7} with $\mu$ not rotationally invariant
   and define $(a_n)_{n=0}^{+\infty}$ as in \eqref{2.12.5.7}. Thus
   there is and $s\in\mathbb R$ such that $\mu(\tau)\neq\mu(\E
   ^{\I s}\tau)$ for some subset $\tau$ of $\mathbb C$; make
   $\tau$ maximal closed with respect to this property. Let $\nu$
   be a measure on $\mathbb C$ transported from $\mu$ via the
   rotation $z\to\E^{-\I s}z$ and let $m_\nu$ be the the measure
   on $[0,+\infty)$ constructed from $\nu$ in the way $m_\mu$ was
   from $\mu$, cf. \eqref{2.12.5.7}. Because, what is a matter of
   straightforward calculation, $m_\mu$ and $m_\nu$ differ on
   $\zb{|z|}{z\in\tau}$, we get indeterminacy of
   $(a_n)_{n=0}^{+\infty}$ at once.

   The other way around, if $m_1$ and $m_2$ are two different
   measures on $[0,+\infty)$ representing the Stieltjes moment
   sequence $(a_n)_{n=0}^{+\infty}$ satisfying \eqref{3.12.5.7},
   then the measure $\mu$ on $\mathbb C$ defined by
   \begin{multline*}
   \mu(\sigma)\okr(2\pi)^{-1}(s\int_0^{a}\D
   t\int_0^{+\infty}\chi_\sigma(r\E^{\I t})m_1(\D
   r)+(1-s)\int_a^{2\pi}\D t\int_0^{+\infty}\chi_\sigma(r\E^{\I
   t})(sm_2(\D r),\\ \sigma \text{ Borel subset of } \mathbb
   C,\;0<s<1,\;0<a<2\pi
   \end{multline*}
   is not rotationally invariant while still \eqref{11.14.5.7} is
   maintained.
 \end{proof}

   \subsection*{R\'esum\'e.} Define two linear operators
$M$ and $D_q$ acting on functions
   \begin{equation}\label{2.6.5.7}
   (Mf)(z)\okr zf(z),\quad (D_qf)(z)\okr
   \begin{cases}
   \frac{f(z)-f(qz)}{z-qz} & \text{if $q\neq1$} \\
   f^\prime(z) & \text{if $q=1$}.
   \end{cases}
   \end{equation}
   It turns out that for $a_+=M$ and $a_-=D_q$ the commutation
   relation \eqref{w28.2} is always satisfied. What Bargmann did
   in \cite{barg} was to find, for $q=1$, a Hilbert space of
   entire functions such that $M$ and $D_1$ are formally adjoint.
   This for arbitrary $q>0$ leads to the reproducing kernel
   Hilbert space $\hhc_q$ of analytic functions with the kernel
   \begin{equation*}
    K(z,w)\okr e_q((1-q)z\bar w) \quad z,w\in|1-q|^{-1/2}\omega_q
   \end{equation*}
   where
   \begin{equation*}
   \omega_q=
   \begin{cases}
   \zb{z}{|z|<1} & \text{if $0<q<1$} \\
   \mathbb C& \text{if $q>1$}
   \end{cases}
   \end{equation*}
    Under these circumstances we always have
   \begin{equation*}
   \is{Z^m}{Z^n}_{\hhc_q}=\delta_{m,n}[m]_q!
   \end{equation*}
   and the operator $S=M$ act as a weighted shift with the weights
   $(\sqrt{[n+1]_q})$ as in Sample Theorem on p. \pageref{b.1}.

Our keynote, subnormality of $M$ now means precisely
\eqref{11.14.5.7} with some $\mu$ is retained. Here we have three
qualitatively different situations:
   \begin{enumerate}
   \item[(a)] for $0<q< 1$ the multiplication operator $M$  is
bounded and subnormal, this implies uniqueness of $\mu$;
   \item[(b)] for $q=1$ the multiplication operator is unbounded
and subnormal, it has a normal extension of cyclic type in the
sense of \cite{try3} and consequently $\mu$ is uniquely determined
as well;
   \item[(c)] for $q>1$   the multiplication operator is unbounded
and subnormal, it has no normal extension of cyclic type in the
sense of \cite{try3} though it does plenty of those of spectral
type in the sense of \cite{try3}, which are not unitary equivalent
\,\footnote{\;That is, there is no unitary map between the
$\mathcal L^2$ spaces in question, which is the identity on
$\hhc_q$.}; explicit example of such, based on \cite{askey}, can
be found in \cite{exp} (one has to replace $q$ by $q\od$ there to
get the commutation relation \eqref{w28.2} satisfied), an explicit
example of non radially invariant measure $\mu$ is struck out in
\cite{kro} and it also comes out from Theorem \ref{t3.12.5.7}.
   \end{enumerate}
   \vspace{1cm}
   \subsection*{The author's afterword.} The fundamentals of this
paper have been presented on several occasions for the last couple
of years, recently at the B\c{e}dlewo 9th Workshop {\it
Noncommutative Harmonic Analysis with Applications to
Probability}. It was Marek Bo\.zejko's contagious enthusiasm what
catalysed converting at long last my distracted notes into a
cohesive exposition.


    \bibliographystyle{amsplain}
    
    \end{document}